# Sur l'accélération de la convergence de la « série de Mādhava-Leibniz »

David Pouvreau[1]

## Résumé

Cet article présente des résultats très novateurs obtenus entre le milieu du XIVe siècle et le début du XVIe siècle par des astronomes indiens de l'école dite « de Mādhava ». Ces résultats, qui s'inscrivent dans le cadre de recherches trigonométriques, concernent la rectification du huitième de circonférence d'un cercle. Ils exposent non seulement un analogue du développement en série de arctan(1), en général connu sous le nom de « série de Leibniz », mais aussi d'autres analogues de développements en série dont la convergence est beaucoup plus rapide. Ces développements sont dérivés d'évaluations des restes des sommes partielles de la série initiale au moyen de réduites de fractions continues généralisées. Une justification en termes modernes en est fournie, qui vise à en restituer tout l'intérêt mathématique.

## Abstract

This paper expounds very innovative results achieved between the mid-14th century and the beginning of the 16th century by Indian astronomers belonging to the so-called "Mādhava school". These results were in keeping with researches in trigonometry: they concern the calculation of the eight of the circumference of a circle. They not only expose an analog of the series expansion of arctan(1) usually known as the "Leibniz series", but also other analogs of series expansions, the convergence of which is much faster. These series expansions are derived from evaluations of the rests of the partial sums of the primordial series, by means of some convergents of generalized continued fractions. A justification of these results in modern terms is provided, which aims at restoring their full mathematical interest.

## Introduction

Bien que le rôle historique joué à l'époque médiévale par les Indiens en matière de système de numération soit reconnu, l'histoire des mathématiques indiennes demeure dans son ensemble largement ignorée. Ce fait a de multiples raisons, parmi lesquelles l'obstacle linguistique. Il est pourtant au moins un autre domaine où ces mathématiques ont aussi eu une importance significative : celui de la trigonométrie. Cette science naquit en Inde vers le IVe siècle après J.-C., vraisemblablement à la suite de l'introduction d'éléments d'astronomie ptolémaïque, consécutive aux contacts établis entre les mondes indien et hellénistique. Elle fut considérablement développée en Inde au cours des douze siècles suivants. L'importance de cette science pour les Indiens avait des raisons primordialement religieuses, comme l'illustre la combinaison des trois citations suivantes :

> « Les Vedas ont été révélés afin d'assurer l'accomplissement des sacrifices ; mais ces sacrifices ont été instaurés en fonction de périodes précises. Pour cette raison, seul celui qui connaît l'astronomie, la science du temps, comprend les sacrifices. » (Rig-Vedāṅga Jyotiṣa, IVe siècle avant J.-C.) [**10**, 1].

[1] Professeur agrégé de mathématiques au Lycée Gustave Eiffel de Bordeaux. Docteur en histoire des sciences, chercheur associé au département de philosophie de l'Université Toulouse II. Email : david_pouvreau@orange.fr

> « Une bonne connaissance de cette science du calcul comprenant la détermination de la position des planètes, l'arithmétique et l'algèbre, constitue les racines de l'arbre de la science des astres. » (Mahādeva, XIII^e siècle) [**10**, 2].

> « Le titre de maître en astronomie sera décerné à celui qui a acquis une connaissance suffisante de la trigonométrie. » (Bhāskara II, XII^e siècle) [**1**, 230].

Parmi les nombreux apports indiens à la trigonométrie, mentionnons l'introduction des demi-cordes (au lieu des cordes utilisées par les Grecs), l'énoncé de toutes les formules trigonométriques fondamentales et la construction de tables très précises (dont certaines furent reprises par les Arabes). Les termes latins sinus et cosinus ont d'ailleurs probablement pour origine les termes sanscrits correspondants (*jyā* et *koṭijyā*). Je me propose toutefois ici d'exposer brièvement quelques résultats obtenus en Inde médiévale qui présentent un intérêt particulier : leur portée dépasse en effet largement le domaine de la trigonométrie en tant que telle pour concerner directement ce qui allait devenir en Europe quelques siècles plus tard la théorie des fonctions analytiques.

Entre 1340 et 1425 vécut au Kerala (une région du Sud-Ouest de l'Inde) un astronome nommé Mādhava de Sangamagrama, désormais considéré comme l'un des plus brillants mathématiciens indiens. Ses traités ont pour l'essentiel disparu et son œuvre n'est connue que par l'intermédiaire de certains astronomes qui se réclamèrent de son enseignement. On désigne par école de Mādhava la lignée de maîtres et de disciples qui se succédèrent du début du XV^e à la fin du XVI^e siècle afin de préserver et d'enrichir cet enseignement. Les membres les plus renommés de cette école sont Parameśvara (ayant vécu entre 1360 et 1460), Nīlakaṇṭha (né en 1444), Jyeṣṭhadeva (né vers 1500) et Śāṇkara Vāriyar (né vers 1500). Leurs œuvres présentent un intérêt considérable, tant du point de vue de l'astronomie que des mathématiques, parce qu'elles révèlent une approche très originale de certains problèmes.

Il en est ainsi de leur calcul du rapport de la circonférence d'un cercle à son diamètre et de leur construction des tables trigonométriques. L'école de Mādhava fut en effet la première dans l'histoire, et ce près de trois siècles avant les Européens, à effectuer une approche de ces problèmes par des « développements en série ». Nīlakaṇṭha, Jyeṣṭhadeva et Śāṇkara Vāriyar énoncèrent et justifièrent ainsi des analogues des développements en série entière de l'arctangente, du sinus et du cosinus. Ils attribuèrent ces énoncés à Mādhava, ainsi que les résultats numériques qui en furent déduits, à savoir d'excellentes tables trigonométriques (correspondant à des valeurs du sinus correctes à $10^{-6}$ près) et l'approximation $\frac{2827433388233}{900000000000}$ du nombre $\pi$, qui est correcte à $2{,}5 \times 10^{-12}$ près.

Je n'exposerai pas ici la manière aussi ingénieuse que complexe dont ces résultats ont été obtenus, ce qui ne saurait se faire en quelques pages. Le lecteur intéressé par le détail des raisonnements et techniques utilisés est invité à consulter le livre que j'ai publié en 2003 à ce sujet [**7**]. L'objectif de cet article se limitera à retranscrire les énoncés indiens relatifs à l'analogue du développement en série de arctan(1) et à l'analyse qui en fut faite, à savoir la mise en œuvre d'analogues de techniques d'accélération de convergence de séries numériques, en vue d'obtenir des approximations de $\pi$ telles que celle énoncée plus haut. Mais il s'agit aussi bien, par-delà cet aspect historique, de considérer de quelle manière ces résultats le plus souvent énoncés sans justification à leur époque peuvent être démontrés et généralisés en termes modernes. Il s'agira donc aussi ici de mettre ainsi en évidence l'intérêt purement mathématique de ces résultats, même du point de vue contemporain.

## 1 – La « série de Mādhava-Leibniz »

Dans le cadre de ses recherches sur la quadrature du cercle, Gottfried Wilhelm Leibniz énonça en 1673 le développement en série que l'on écrit désormais :

$$\frac{\pi}{4} = \sum_{p=1}^{+\infty} \frac{(-1)^{p-1}}{2p-1} \qquad (1)$$

Le mathématicien allemand semble avoir été le premier Européen à l'énoncer, bien que James Gregory ait établi deux ans plus tôt la formule plus générale donnant le développement en série entière de l'arctangente. Il est bien connu que ce résultat peut par exemple s'obtenir comme suit :

$$\forall\, n \geq 1, \qquad \frac{\pi}{4} = \arctan(1) = \int_0^1 \frac{dx}{1+x^2} = \int_0^1 \left( \sum_{p=0}^{n-1} (-1)^p x^{2p} \right) dx + (-1)^n \int_0^1 \frac{x^{2n}}{1+x^2} dx$$

$$= \sum_{p=1}^{n} \frac{(-1)^{p-1}}{2p-1} + (-1)^n \rho_n$$

$$\text{avec } 0 < \rho_n = \int_0^1 \frac{x^{2n}}{1+x^2} dx \leq \int_0^1 x^{2n} dx = \frac{1}{2n+1} \qquad \text{et par conséquent : } \lim_{n \to +\infty} \rho_n = 0$$

En réalité, plusieurs disciples de Mādhava avaient déjà donné sous diverses formes des évaluations du huitième de la circonférence d'un cercle à partir d'un « développement en série » analogue à (1). Dans le *Yuktidīpikā* (littéralement : « éclairage du fondement »), Śāṇkara Vāriyar énonça par exemple, en attribuant cette règle à Mādhava lui-même :

> « Prends le diamètre du cercle multiplié par 4 et divisé par 1 ; soustrais et ajoute alternativement à ce résultat les termes consécutifs obtenus en divisant le quadruple du diamètre plusieurs fois par les nombres impairs 3, 5, etc…[…] Le résultat est la circonférence du cercle. En prenant plus de termes, le résultat sera plus précis. » (Śāṇkara Vāriyar, milieu du XVI[e] siècle) [**8**, 94], [**6**, 68-69].

Si l'on note *D* le diamètre du cercle de circonférence *C*, cet énoncé peut être retranscrit par :

$$C \simeq \frac{4D}{1} - \frac{4D}{3} + \frac{4D}{5} - \frac{4D}{7} + \cdots + \frac{(-1)^{n-1} 4D}{2n-1}$$

(cette approximation étant d'autant meilleure que *n* est grand)

Compte tenu de $C = \pi D$, on constate que ce résultat est analogue à (1), qu'il est donc légitime d'appeler le développement en « série de Mādhava-Leibniz » [**6**]. Il faut néanmoins remarquer que l'on peut seulement parler d'analogie, parce qu'il n'est pas question chez les Indiens d'une sommation infinie (actuelle), mais d'une approximation dont la précision est aussi grande que voulue.

L'approximation de $\pi$ attribuée à Mādhava par ses disciples énoncée plus haut est, au moins pour l'époque médiévale, d'une qualité exceptionnelle. Elle fut obtenue à partir du développement en série de Mādhava-Leibniz. Toutefois, l'utilisation telle quelle de cette série ne saurait expliquer la qualité de cette approximation. Sa convergence est en effet excessivement lente : il faut sommer plus de 2000 termes afin d'obtenir ne serait-ce que trois décimales significatives. En réalité, et là est tout l'intérêt de ce qui va suivre, les mathématiciens de l'école de Mādhava ont trouvé divers procédés permettant d'accélérer la convergence de cette série, qu'il s'agit ici d'exposer.

## 2 – L'évaluation de la valeur absolue des restes des sommes partielles de la série de Mādhava-Leibniz

Pour tout $n \geq 1$, notons :

$$S_n = \sum_{p=1}^{n} \frac{(-1)^{p-1}}{2p-1}$$

Il a précédemment été montré que :

$$\forall\, n \geq 1, \qquad \frac{\pi}{4} = S_n + (-1)^n \rho_n \qquad \text{avec} \quad \rho_n = \int_0^1 \frac{x^{2n}}{1+x^2}\,dx$$

Une bonne évaluation $R_n$ de $\rho_n$ (valeur absolue du reste de la somme partielle de la série au rang *n*) peut donc permettre d'accélérer la convergence, en ajoutant un terme correcteur aux sommes partielles.

## 2-1 – Les évaluations de la valeur absolue des restes énoncées dans l'école de Mādhava

Examinons d'abord les termes correcteurs déterminés par les membres de l'école de Mādhava. À la suite de l'énoncé figurant plus haut, Śāṇkara Vāriyar expliqua comment une meilleure approximation du huitième de circonférence peut être obtenue :

> « […] Prends le nombre pair immédiatement supérieur au nombre impair auquel le procédé précédent a été interrompu. Comme auparavant, multiplie le quadruple du diamètre par la moitié de ce nombre pair et divise par son carré ajouté de 1. Le quotient doit être additionné à la série si le dernier terme a été soustrait et soustrait si le dernier terme a été additionné. » (Śāṇkara Vāriyar) [**8**, 94], [**6**, 68-69].

C'est-à-dire, en conservant les notations introduites plus haut :

$$\text{Si } n \text{ est grand, alors : } \; C \simeq \frac{4D}{1} - \frac{4D}{3} + \frac{4D}{5} - \frac{4D}{7} + \cdots \pm \frac{4D}{2n-1} \mp \frac{4D \times \left(\frac{2n}{2}\right)}{(2n)^2+1}$$

Soit encore :

$$\text{Si } n \text{ est grand, alors : } \; \frac{\pi}{4} \simeq 1 - \frac{1}{3} + \frac{1}{5} - \frac{1}{7} + \cdots \pm \frac{1}{2n-1} \mp R_n^{(2)} \qquad \text{avec} \quad R_n^{(2)} = \frac{n}{4n^2+1}$$

(cette approximation étant d'autant meilleure que *n* est grand)

Śāṇkara Vāriyar énonça ensuite :

> « Pour une plus grande précision […], dans le terme final, le multiplicande du quadruple du diamètre est le carré de la moitié du nombre pair additionné à 1 et le diviseur est le quadruple de ce multiplicande additionné à 1, multiplié ensuite par la moitié du nombre pair. » (Śāṇkara Vāriyar) [**8**, 94], [**6**, 68-69].

C'est-à-dire :

$$\text{Si } n \text{ est grand, alors : } \; C \simeq \frac{4D}{1} - \frac{4D}{3} + \frac{4D}{5} - \frac{4D}{7} + \cdots \pm \frac{4D}{2n-1} \mp \frac{4D \times \left[\left(\frac{2n}{2}\right)^2 + 1\right]}{\left[4\left(\left(\frac{2n}{2}\right)^2 + 1\right) + 1\right] n}$$

Ou encore :

$$\text{Si } n \text{ est grand, alors : } \; \frac{\pi}{4} \simeq 1 - \frac{1}{3} + \frac{1}{5} - \frac{1}{7} + \cdots \pm \frac{1}{2n-1} \mp R_n^{(3)} \qquad \text{avec} \quad R_n^{(3)} = \frac{n^2+1}{(4n^2+5)n}$$

(cette approximation étant meilleure que la précédente)

Śāṇkara Vāriyar, au cours des explications qu'il fournit quant à la manière d'obtenir ces termes correcteurs, donna aussi $R_n^{(1)} = 1/4n$ pour première approximation des restes [**5**, 150].

La qualité de ces évaluations est remarquable. Il suffit pour s'en rendre compte de l'inscrire dans le cadre de la théorie des développements en fractions continues.

## 2-2 – Justification et généralisation des évaluations énoncées dans l'école de Mādhava

Rappelons qu'une fraction continue est la donnée de deux suites $(a_n)$ et $(b_n)$ et de l'expression convergente

$$\cfrac{a_0}{b_0+\cfrac{a_1}{b_1+\cfrac{a_2}{b_2+\cfrac{a_3}{b_3+\cdots}}}}$$

Si $a_n = 1$ pour tout $n \geq 1$, on dit qu'il s'agit d'une fraction continue simple. Sinon, on parle de fraction continue généralisée. $a_0/b_0$ est la « première réduite » de la fraction continue, $\frac{a_0}{b_0+a_1/b_1}$ la seconde réduite, etc. Si $a_n > 0$ et $b_n > 0$ pour tout $n \geq 1$, alors les réduites sont des approximations alternativement par excès et par défaut de la limite.

Revenons maintenant à l'évaluation de la valeur absolue $\rho_n$ du reste d'une somme partielle de la série de Mādhava-Leibniz. On peut montrer que cette valeur absolue a pour développement en fraction continue généralisée :

$$\rho_n = \cfrac{1/2}{2n+\cfrac{1^2}{2n+\cfrac{2^2}{2n+\cfrac{3^2}{2n+\cdots}}}}$$

Et l'on constate dès lors que les trois évaluations de $\rho_n$ données par Śāṇkara Vāriyar, notées ici $R_n^{(1)}, R_n^{(2)}$ et $R_n^{(3)}$, sont en fait respectivement la première, la seconde et la troisième réduite de cette fraction continue généralisée. En effet :

$$R_n^{(1)} = \frac{1/2}{2n} \quad ; \quad R_n^{(2)} = \frac{1}{4n+\frac{1}{n}} = \frac{1/2}{2n+\frac{1^2}{2n}} \quad ; \quad R_n^{(3)} = \cfrac{1}{4n+\cfrac{1}{n+\frac{1}{n}}} = \cfrac{1/2}{2n+\cfrac{1^2}{2n+\frac{2^2}{2n}}}$$

La manière dont ces évaluations de la valeur absolue des restes ont été déterminées par les Indiens ne sera pas examinée ici. Il suffira de dire qu'elles furent essentiellement obtenues par induction à partir d'approximations du rapport de la circonférence au diamètre déjà connues (notamment $\frac{62832}{20000}$ ) et qu'une reconstruction rationnelle *a posteriori* fut tentée par Śāṇkara Vāriyar [**7**, 72-76]. La procédure alors utilisée tint très probablement à l'observation du fait qui s'exprime en termes modernes par :

$$\forall\, n \in [\![1;4]\!], \exists\, f_n \in\, ]0;1[, \left|\frac{1}{4}\times\frac{62832}{20000} - \sum_{p=1}^{n}\frac{(-1)^{p-1}}{2p-1}\right| = \cfrac{1}{4n+\cfrac{1}{n+\cfrac{1}{n+f_n}}}$$

Le développement en fraction continue généralisée annoncé peut être obtenu comme suit. Remarquons d'abord que, pour tout $n \geq 1$ :

$$\rho_n + \rho_{n+1} = (-1)^n\left[\left(\frac{\pi}{4} - S_n\right) - \left(\frac{\pi}{4} - S_{n+1}\right)\right] = (-1)^n(S_{n+1} - S_n) = (-1)^n\frac{(-1)^n}{2n+1}$$

Par conséquent :

$$\forall\, n \geq 1, \qquad \rho_n + \rho_{n+1} = \frac{1}{2n+1} \qquad (2)$$

$(\rho_n)_{n\geq 1}$ est l'unique suite convergente vérifiant (2).
En effet, soit $(r_n)_{n\geq 1}$ une suite vérifiant (2). Posons $v_n = r_n - \rho_n$ pour tout $n \geq 1$. Il est clair que la suite de terme général $(v_{n+1} + v_n)$ est constamment nulle. Il en résulte : $\forall\, n \geq 1,\ v_n = (-1)^{n-1}v_1$. C'est-à-dire aussi : $\forall\, n \geq 1,\ r_n = \rho_n + (-1)^{n-1}v_1$. Par conséquent :

$$(r_n) \text{ converge} \Leftrightarrow \lim_{n\to+\infty} (-1)^{n-1}v_1 \text{ existe et est réelle} \Leftrightarrow v_1 = 0 \Leftrightarrow \forall\, n \geq 1,\ r_n = \rho_n$$

Supposons maintenant qu'il existe une suite $(a_k)_{k\in\mathbb{N}}$ telle que, pour tout $n \geq 1$ :

$$\rho_n = \cfrac{a_0}{n + \cfrac{a_1}{n + \cfrac{a_2}{n + \cdots}}}$$

Pour tout $k \in \mathbb{N}$, notons alors :

$$F_n^{(k)} = \cfrac{a_k}{n + \cfrac{a_{k+1}}{n + \cfrac{a_{k+2}}{n + \cdots}}}$$

Pour tout $k \in \mathbb{N}$, $a_k/n$ est un équivalent de $F_n^{(k)}$ lorsque *n* est voisin de $+\infty$. Par conséquent :

$$\lim_{n\to+\infty} F_n^{(k)} = 0$$

Par ailleurs, il résulte de l'identité (2) que la suite $(a_k)_{k\in\mathbb{N}}$ doit nécessairement vérifier :

$$\forall\, n \geq 1, \qquad \frac{a_0}{n + F_n^{(1)}} + \frac{a_0}{n + 1 + F_{n+1}^{(1)}} = \frac{1}{2n+1}$$

Après réduction au même dénominateur, on obtient :

$$\left(n + F_n^{(1)}\right)\left(n + 1 + F_{n+1}^{(1)}\right) = a_0(2n+1)(2n+1+F_n^{(1)} + F_{n+1}^{(1)})$$

D'où aussi, en ordonnant et en divisant par *n* la relation obtenue :

$$(1-4a_0)n + \left[(1-4a_0) + (1-2a_0)(F_n^{(1)} + F_{n+1}^{(1)})\right] + \frac{1}{n}\left[(1-a_0)F_n^{(1)} - a_0F_{n+1}^{(1)} + F_n^{(1)}F_{n+1}^{(1)} - a_0\right] = 0$$

Puisque $\lim_{n\to+\infty} F_n^{(1)} = 0 = \lim_{n\to+\infty} F_{n+1}^{(1)}$, il est donc nécessaire que $1 - 4a_0 = 0$, soit $a_0 = 1/4$.
L'identité précédente devient alors :

$$2n\left(F_n^{(1)} + F_{n+1}^{(1)}\right) + 3F_n^{(1)} - F_{n+1}^{(1)} + 4F_n^{(1)}F_{n+1}^{(1)} - 1 = 0 \qquad (3)$$

Or, pour tout $n \geq 1$ :

$$F_n^{(1)} = \frac{a_1}{n + F_n^{(2)}} \qquad \text{et} \qquad F_{n+1}^{(1)} = \frac{a_1}{n + 1 + F_{n+1}^{(2)}}$$

En reportant dans (3), on obtient donc, après réduction au même dénominateur :

$$2na_1\left(n + F_n^{(2)}\right) + 2na_1\left(n + 1 + F_{n+1}^{(2)}\right) - \left(n + F_n^{(2)}\right)\left(n + 1 + F_{n+1}^{(2)}\right) + 3a_1\left(n + 1 + F_{n+1}^{(2)}\right) - a_1\left(n + F_n^{(2)}\right) + 4{a_1}^2 = 0$$

D'où, après avoir ordonné et divisé par *n* la relation obtenue :

$$(1-4a_1)n+\left[(1-4a_1)+(1-2a_1)(F_n^{(2)}+F_{n+1}^{(2)})\right]$$
$$+\frac{1}{n}\left[(1+a_1)F_n^{(2)}-3a_1F_{n+1}^{(2)}+F_n^{(2)}F_{n+1}^{(2)}-4{a_1}^2-3a_1\right]=0$$

Il résulte alors de $\lim_{n\to+\infty} F_n^{(2)}=0=\lim_{n\to+\infty} F_{n+1}^{(2)}$ que $1-4a_1=0$ nécessairement, soit : $a_1=1/4$.
En reportant cette valeur dans l'égalité précédente, on obtient :

$$2n\left(F_n^{(2)}+F_{n+1}^{(2)}\right)+5F_n^{(2)}-3F_{n+1}^{(2)}+4F_n^{(2)}F_{n+1}^{(2)}-4=0$$

Pour tout entier $k\geq 2$, notons maintenant $(H_k)$ l'hypothèse de récurrence :

$$\begin{cases} a_{k-1}=\dfrac{(k-1)^2}{4} \\ \forall\, n\geq 1, \qquad 2n\left(F_n^{(k)}+F_{n+1}^{(k)}\right)+(2k+1)F_n^{(k)}-(2k-1)F_{n+1}^{(k)}+4F_n^{(k)}F_{n+1}^{(k)}-k^2=0 \end{cases}$$

Il a été établi plus haut que $(H_2)$ est vraie. Supposons $(H_k)$ vraie pour un certain $k\geq 2$.
En utilisant

$$F_n^{(k)}=\frac{a_k}{n+F_n^{(k+1)}} \quad \text{et} \quad F_{n+1}^{(k)}=\frac{a_k}{n+1+F_{n+1}^{(k+1)}}$$

on obtient :

$$2na_k\left(n+F_n^{(k+1)}\right)+2na_k\left(n+1+F_{n+1}^{(k+1)}\right)-k^2\left(n+F_n^{(k+1)}\right)\left(n+1+F_{n+1}^{(k+1)}\right)$$
$$+(2k+1)a_k\left(n+1+F_{n+1}^{(k+1)}\right)-(2k-1)a_k\left(n+F_n^{(k+1)}\right)+4{a_k}^2=0$$

D'où, après avoir développé, ordonné et divisé par *n* la relation obtenue :

$$(k^2-4a_k)n+\left[(k^2-4a_k)+(k^2-2a_k)(F_n^{(k+1)}+F_{n+1}^{(k+1)})\right]$$
$$+\frac{1}{n}\Big[(k^2+(2k-1)a_k)F_n^{(k+1)}-(2k+1)a_kF_{n+1}^{(k+1)}+k^2F_n^{(k+1)}F_{n+1}^{(k+1)}-4{a_k}^2$$
$$-(2k+1)a_k\Big]=0$$

Il résulte alors de $\lim_{n\to+\infty} F_n^{(k+1)}=0=\lim_{n\to+\infty} F_{n+1}^{(k+1)}$ que $k^2-4a_k=0$ nécessairement, soit $a_k=k^2/4$.
L'identité précédente devient alors :

$$2n\left(F_n^{(k+1)}+F_{n+1}^{(k+1)}\right)+(2k+3)F_n^{(k+1)}-(2k+1)F_{n+1}^{(k+1)}+4F_n^{(k+1)}F_{n+1}^{(k+1)}-(k+1)^2=0$$

Ceci établit que $(H_{k+1})$ est vraie. Par récurrence, $(H_k)$ est donc vraie pour tout $k\geq 2$. D'où résulte que nécessairement : $a_k=k^2/4$ pour tout $k\geq 1$.

Si $\rho_n$ admet un développement en fraction continue du type proposé, il ne peut donc être que

$$\rho_n=\cfrac{1/4}{n+\cfrac{1^2/4}{n+\cfrac{2^2/4}{n+\cfrac{3^2/4}{n+\cdots}}}}$$

Ce qui correspond au développement annoncé. Réciproquement, le développement en fraction continue précédent vérifie par construction la relation (2). Comme il converge (vers 0) lorsque *n* tend vers l'infini, le fait que $(\rho_n)_{n\geq 1}$ est l'unique suite convergente vérifiant (2) montre que ce développement est bien celui de $\rho_n$. A de la sorte été démontré le résultat suivant :

**Proposition 1.**

$$\forall\, n \geq 1, \qquad \pi = 4\sum_{p=1}^{n} \frac{(-1)^{p-1}}{2p-1} + (-1)^n \cfrac{2}{2n + \cfrac{1^2}{2n + \cfrac{2^2}{2n + \cfrac{3^2}{2n + \cdots}}}}$$

Remarquons, avec $n = 1$, le cas particulier :

$$\pi = 4 - \cfrac{2}{2 + \cfrac{1^2}{2 + \cfrac{2^2}{2 + \cfrac{3^2}{2 + \cdots}}}}$$

On peut observer la similitude de ce résultat avec une formule restée fameuse [**2**, p. 10] qui fut énoncée sans justification par William Brouncker au moins un siècle plus tard, en 1655.

Comme l'indiquent les valeurs numériques présentées plus loin, la correction des sommes partielles de la série de Mādhava-Leibniz par des réduites de $\rho_n$ est très performante. Il faut remarquer qu'une telle efficacité est *a priori* étonnante, car aussi bien la série en tant que telle que le développement en fraction continue généralisée en tant que tel ont, *séparément*, une convergence très lente.

## 3 – Une méthode de construction de nouvelles séries à convergence accélérée

Je vais maintenant considérer de quelle manière de nouvelles séries convergeant beaucoup plus rapidement que la série initiale peuvent être construites à partir des termes correcteurs ainsi déterminés.

### 3-1 – Les développements à convergence accélérée énoncés dans l'école de Mādhava

Śāṇkara Vāriyar donna dans le *Yuktidīpikā* d'autres « développements en série » de la circonférence que ceux exposés plus haut, en les attribuant là encore à Mādhava. Ces développements convergent beaucoup plus vite que la série initiale. Śāṇkara Vāriyar ne les justifia pas (aucun auteur indien de l'époque ne le fit d'ailleurs), mais ils furent donnés à la suite des énoncés concernant les termes correcteurs retranscrits au paragraphe 2 et en résultent donc manifestement :

> « La circonférence est de façon analogue obtenue lorsque quatre fois le diamètre est divisé par les cubes des nombres impairs en commençant par 3 diminués par ces nombres eux-mêmes et que ces quotients sont alternativement additionnés ou soustraits au triple du diamètre. » (Śāṇkara Vāriyar) [**8**, 95], [**1**, 266].

C'est-à-dire, en conservant les notations introduites plus haut :

$$C = 3D + \frac{4D}{3^3 - 3} - \frac{4D}{5^3 - 5} + \frac{4D}{7^3 - 7} - \frac{4D}{9^3 - 9} + \cdots$$

Ce qui peut être interprété en termes modernes par l'identité :

$$\pi = 3 + 4\sum_{p=1}^{+\infty} \frac{(-1)^{p-1}}{(2p+1)^3 - (2p+1)}$$

Le *Yuktidīpikā* énonça aussi :

> « Le quadruple des nombres impairs est ajouté à leur cinquième puissance ; 16 fois le diamètre est successivement divisé par les nombres obtenus ; les résultats des rangs impairs sont additionnés et ceux de rang pair soustraits. La circonférence correspondant au diamètre est ainsi obtenue. » (Śāṇkara Vāriyar) [**8**, 95].

C'est-à-dire :

$$C = \frac{16D}{1^5 + 4 \times 1} - \frac{16D}{3^5 + 4 \times 3} + \frac{16D}{5^5 + 4 \times 5} - \frac{16D}{7^5 + 4 \times 7} + \cdots$$

À cet énoncé correspond clairement le développement en série :

$$\pi = 16 \sum_{p=0}^{+\infty} \frac{(-1)^p}{(2p+1)^5 + 4(2p+1)}$$

Le problème consiste donc à savoir comment ces développements ont été obtenus par les astronomes kéralais. Je me propose ici de montrer une propriété générale explicitant la manière dont on peut les déduire de l'évaluation asymptotique des restes des sommes partielles de la série de Mādhava-Leibniz : seront ainsi en particulier justifiés les deux « développements en série » énoncés par Śāṇkara Vāriyar.

### 3-2 – Justification des développements à convergence accélérée énoncés dans l'école de Mādhava : une méthode générale de construction

$(S_n)$ étant une suite réelle convergeant vers un réel $\varphi$, on dit qu'une suite réelle $(S'_n)$ converge vers $\varphi$ plus rapidement que $(S_n)$ si :

$$\lim_{n \to +\infty} \frac{\varphi - S'_n}{\varphi - S_n} = 0$$

Une suite réelle $(u_n)_{n\geq 1}$ étant donnée telle que la série de somme partielle

$$S_n = \sum_{p=1}^{n} (-1)^{p-1} u_p$$

converge vers $\varphi$, on cherche dans le cas présent à construire une série de somme partielle $S'_n$ convergeant vers $\varphi$ plus rapidement que $(S_n)$. Étudions d'abord « constructivement » ce problème.

Posons $S'_n = S_n + (-1)^n R_n$, où $R_n$ est un équivalent de $\rho_n = (-1)^n(\varphi - S_n)$. Supposons alors qu'il existe une suite $(v_n)_{n\geq 1}$ telle que pour tout $n \geq 1$ :

$$S'_n = \sum_{p=1}^{n} (-1)^{p-1} v_p$$

Pour tout $n \geq 1$, on a alors d'une part : $S'_{n+1} - S'_n = (-1)^n v_{n+1}$ ; et d'autre part :

$$S'_{n+1} - S'_n = (S_{n+1} - S_n) + (-1)^{n+1} R_{n+1} - (-1)^n R_n = (-1)^n [u_{n+1} - (R_{n+1} + R_n)]$$

On a donc nécessairement : $v_n = u_n - (R_n + R_{n-1})$ pour tout $n \geq 2$.

Par ailleurs, $S'_1 = v_1 = S_1 - R_1 = u_1 - R_1$. Il en résulte finalement :

$$S'_n = u_1 - R_1 + \sum_{p=2}^{n} (-1)^{p-1} [u_p - (R_p + R_{p-1})] = u_1 - R_1 + \sum_{q=1}^{n-1} (-1)^{q-1} [(R_{q+1} + R_q) - u_{q+1}]$$

Pour chaque $n \geq 1$, notons par ailleurs $\varepsilon_n$ le nombre $\left(\frac{R_n}{\rho_n} - 1\right)$, qui tend vers 0. La suite $(\varepsilon_n)_{n\geq 1}$ sera appelée ici la « qualité » de la suite d'équivalents $(R_n)_{n\geq 1}$ . Une « qualité » sera alors dite « meilleure » qu'une autre si elle converge plus vite vers 0. On peut dans ces conditions énoncer le résultat suivant :

**Proposition 2.** Soit $(u_p)_{p\geq 1}$ une suite réelle telle que la série de terme général $(-1)^{p-1}u_p$ converge vers un réel $\varphi$. On note $(S_n)_{n\geq 1}$ la suite des sommes partielles de cette série et $\rho_n = (-1)^n(\varphi - S_n)$ pour tout $n \geq 1$. Soit $(R_n)_{n\geq 1}$ une suite réelle de limite nulle. On note $(S'_n)_{n\geq 1}$ la suite dont le terme général est $S'_n = S_n + (-1)^n R_n$ , $(v_p)_{p\geq 1}$ la suite de terme général $v_p = \left(R_p + R_{p+1}\right) - u_{p+1}$ et $(S_n'')_{n\geq 1}$ la suite de terme général

$$S''_n = u_1 - R_1 + \sum_{p=1}^{n} (-1)^{p-1} v_p$$

(i) Les suites $(S'_n)$ et $(S''_n)$ convergent vers $\varphi$.
(ii) La suite $(S'_n)$ converge vers $\varphi$ plus rapidement que la suite $(S_n)$ si et seulement si $R_n$ est un équivalent de $\rho_n$.
(iii) Si $R_n$ est un équivalent de $\rho_n$, alors la suite $(S''_n)$ converge vers $\varphi$ plus rapidement que $(S_n)$ , et ce avec une vitesse d'autant plus grande que la « qualité » de $(R_n)_{n\geq 1}$ est meilleure.

En effet, il est d'abord clair que $(S'_n)$ converge vers $\varphi$. De plus, pour tout $n \geq 1$ :

$$\begin{aligned} S''_n &= u_1 - R_1 + \sum_{p=1}^{n} (-1)^{p-1}\left(R_p + R_{p+1} - u_{p+1}\right) \\ &= \sum_{p=1}^{n+1} (-1)^{p-1} u_p + \sum_{p=1}^{n} \left(-R_p + R_p\right) + (-1)^{n-1} R_{n+1} \\ &= S_{n+1} + (-1)^{n-1} R_{n+1} = S'_{n+1} \end{aligned}$$

Par conséquent : $\lim_{n\to+\infty} S''_n = \lim_{n\to+\infty} S'_{n+1} = \varphi$. D'où (i). Par ailleurs, observons que pour tout $n \geq 1$ :

$$\frac{\varphi - S'_n}{\varphi - S_n} = \frac{S_n + (-1)^n \rho_n - S_n - (-1)^n R_n}{S_n + (-1)^n \rho_n - S_n} = \frac{\rho_n - R_n}{\rho_n} = 1 - \frac{R_n}{\rho_n}$$

On en déduit que $(S'_n)$ converge vers $\varphi$ plus vite que $(S_n)$ si et seulement si $\lim_{n\to+\infty} R_n/\rho_n = 1$. D'où (ii).
Enfin, supposons que $R_n$ soit un équivalent de $\rho_n$. Notons $(\varepsilon_n)_{n\geq 1}$ la « qualité » de $(R_n)_{n\geq 1}$ , telle que $R_n = (1 + \varepsilon_n)\rho_n$ pour tout $n \geq 1$. On a : $\forall\, n \geq 1, \varphi = S_{n+1} + (-1)^{n+1}\rho_{n+1}$. D'où aussi :

$$\begin{aligned} \forall\, n \geq 1, \qquad |S''_n - \varphi| &= |[S_{n+1} - \varphi] + (-1)^{n-1} R_{n+1}| = |[S_{n+1} - \varphi] + (-1)^{n-1}(1 + \varepsilon_{n+1})\rho_{n+1}| \\ &= |[S_{n+1} - \varphi] + (-1)^{2n}(1 + \varepsilon_{n+1})(\varphi - S_{n+1})| \\ &= |\varepsilon_{n+1}||S_{n+1} - \varphi| \end{aligned}$$

On en déduit que la suite $(S''_n)$ converge vers $\varphi$ plus rapidement que la suite $(S_{n+1})$, et ce avec une vitesse d'autant plus grande que la « qualité » de $(R_n)_{n\geq 1}$ est meilleure. D'où (iii).

Considérons maintenant de ce point de vue la série de Mādhava-Leibniz. On a dans ce cas particulier $u_n = \dfrac{1}{2n-1}$ pour tout $n \geq 1$. L'expression de $\rho_n$ sous forme de fraction continue généralisée a été déterminée au paragraphe 2, ses trois premières réduites ayant été notées $R_n^{(1)}$, $R_n^{(2)}$ et $R_n^{(3)}$.

Avec $R_n^{(1)} = 1/4n$, on trouve sans difficulté que pour tout $n \geq 1$ :

$$v_n = \frac{1}{4n} + \frac{1}{4n+4} - \frac{1}{2n+1} = \frac{1}{(2n+1)^3 - (2n+1)}$$

Or : $u_1 - R_1^{(1)} = 1 - 1/4 = 3/4$. On peut en déduire :

$$\frac{3}{4}+\sum_{p=1}^{+\infty}\frac{(-1)^{p-1}}{(2p+1)^3-(2p+1)}=\frac{\pi}{4}$$

Soit encore :

$$\pi = 3+4\sum_{p=1}^{+\infty}\frac{(-1)^{p-1}}{(2p+1)^3-(2p+1)}$$

Avec $R_n^{(2)} = n/(4n^2+1)$, on obtient que pour tout $n \geq 1$ :

$$v_n = \frac{n}{4n^2+1}+\frac{n+1}{4(n+1)^2+1}-\frac{1}{2n+1}=\frac{4}{(2n+1)^5+4(2n+1)}$$

Or : $u_1 - R_1^{(2)} = 1-1/5 = 4/5$. Par conséquent :

$$\frac{4}{5}-4\sum_{p=1}^{+\infty}\frac{(-1)^{p-1}}{(2p+1)^5+4(2p+1)}=\frac{\pi}{4}$$

Soit encore :

$$\pi = 16\sum_{p=0}^{+\infty}\frac{(-1)^{p}}{(2p+1)^5+4(2p+1)}$$

Il est ainsi rendu compte des deux développements annoncés par Śāṇkara Vāriyar avec des moyens modernes. Le lecteur intéressé trouvera toutefois une justification plus directe et surtout plus conforme aux raisonnements ayant eu cours dans l'école de Mādhava dans le livre déjà mentionné [**7**, 45-79].

Notons que ni Śāṇkara Vāriyar, ni aucun autre membre de cette école, ne semble avoir donné le développement qui résulte de l'utilisation de $R_n^{(3)}$. On obtient dans ce cas par la procédure exposée :

$$\pi = \frac{28}{9}+36\sum_{p=1}^{+\infty}\frac{(-1)^{p-1}}{p(p+1)(2p+1)(4p^2+5)(4p^2+8p+9)}$$

Avant de comparer numériquement la vitesse de convergence de ces trois séries, deux remarques peuvent être faites. La plus évidente est qu'elles satisfont le critère spécial de convergence des séries alternées. De sorte que dans chacun des cas, la majoration de l'erreur est de l'ordre de $1/2n^3$, $1/2n^5$ et $9/8n^7$ respectivement. L'accélération de convergence induite est donc dans une certaine mesure prévisible même sans l'utilisation du (iii) de la propriété établie plus haut.

La seconde remarque porte sur le fait que la première série correspond exactement à celle que l'on obtient après application à la série de Mādhava-Leibniz du très performant et désormais classique algorithme du « delta-2 d'Aitken ». Rappelons que cet algorithme repose sur la propriété suivante :

**Proposition 3 (Delta-2 d'Aitken).** Soit $(S_n)_{n\geq 1}$ une suite convergeant vers un réel $\varphi$. Si la limite du rapport $\frac{\varphi - S_{n+1}}{\varphi - S_n}$ existe et est différente de 1, alors, à condition qu'elle soit bien définie, la suite $(S'_n)_{n\geq 3}$ de terme général $S'_n = \frac{S_n S_{n-2} - {S_{n-1}}^2}{(S_n - S_{n-1})-(S_{n-1}-S_{n-2})}$ converge vers $\varphi$ plus rapidement que $(S_n)_{n\geq 1}$.

Dans le cas présent où $(S_n)_{n\geq 1}$ est la suite des sommes partielles de la série de Mādhava-Leibniz, on obtient pour tout $n \geq 3$ :

$$S_{n-1} = S_n + \frac{(-1)^n}{2n-1} \quad \text{et} \quad S_{n-2} = S_n + (-1)^n\left(\frac{1}{2n-1}-\frac{1}{2n-3}\right)$$

On en déduit par un calcul élémentaire :

$$\forall\, n \geq 3, \qquad S'_n = \frac{S_n S_{n-2} - {S_{n-1}}^2}{(S_n - S_{n-1}) - (S_{n-1} - S_{n-2})} = S_n + (-1)^n \frac{2n-3}{4(n-1)(2n-1)}$$

Notons alors $R_n = \frac{2n-3}{4(n-1)(2n-1)}$ pour tout $n \geq 3$. La propriété démontrée plus haut permet d'affirmer que la suite de terme général

$$S''_n = 1 - \frac{1}{3} + \frac{1}{5} - R_3 + \sum_{p=3}^{n} (-1)^p \left( R_p + R_{p+1} - \frac{1}{2p+1} \right)$$

converge vers $\pi/4$. Or, on montre sans difficulté que pour tout $p \geq 3$ :

$$R_p + R_{p+1} - \frac{1}{2p+1} = -\frac{1}{4p(p-1)(2p-1)}$$

On en déduit que pour tout $n \geq 3$ :

$$S''_n = 1 - \frac{1}{3} + \frac{1}{5} - R_3 + \sum_{p=3}^{n} \frac{(-1)^p}{4p(p-1)(2p-1)} = 1 - \frac{1}{3} + \frac{1}{5} - \frac{3}{40} + \frac{1}{4} \sum_{q=2}^{n-1} \frac{(-1)^{q-1}}{q(q+1)(2q+1)}$$

$$= \frac{19}{24} - \frac{1}{4} \times \frac{1}{6} + \frac{1}{4} \sum_{q=1}^{n-1} \frac{(-1)^{q-1}}{q(q+1)(2q+1)} = \frac{3}{4} + \sum_{q=1}^{n-1} \frac{(-1)^{q-1}}{(2q+1)^3 - (2q+1)}$$

Il en résulte finalement :

$$\frac{3}{4} + \sum_{p=1}^{+\infty} \frac{(-1)^{p-1}}{(2p+1)^3 - (2p+1)} = \frac{\pi}{4}$$

Ceci établit que la série déduite par les astronomes kéralais de l'utilisation de l'équivalent $R_n^{(1)}$ correspond exactement à la série qui se déduit de l'application du « delta-2 d'Aitken ».

## 4- Comparaison numérique et conclusion

Je me limite ici à la comparaison numérique de la convergence des trois séries déterminées précédemment, compte tenu du lien établi entre ces séries et l'addition des termes correcteurs respectifs $R_n^{(1)}$, $R_n^{(2)}$ et $R_n^{(3)}$ aux sommes partielles de la série de Mādhava-Leibniz. Notons, pour tout $n \geq 1$ :

$$a_n = 3 + 4 \sum_{p=1}^{n} \frac{(-1)^{p-1}}{(2p+1)^3 - (2p+1)}$$

$$b_n = 16 \sum_{p=0}^{n} \frac{(-1)^p}{(2p+1)^5 + 4(2p+1)}$$

$$c_n = \frac{28}{9} + 36 \sum_{p=1}^{n} \frac{(-1)^{p-1}}{p(p+1)(2p+1)(4p^2+5)(4p^2+8p+9)}$$

On obtient alors (les caractères gras correspondant aux décimales correctes de $\pi$) :

| | Suite $(a_n)$ | Suite $(b_n)$ | Suite $(c_n)$ |
|---|---|---|---|
| $n$ = 2 | **3,1**333333333333 | **3,1**372549019608 | **3,141**4634146341 |
| $n$ = 3 | **3,14**52380952382 | **3,14**23423423424 | **3,141**6149068323 |
| $n$ = 4 | **3,1**396825396626 | **3,141**391941392 | **3,1415**873015673 |
| $n$ = 5 | **3,14**27128427129 | **3,141**6627377024 | **3,14159**42744802 |
| $n$ = 10 | **3,141**4067184965 | **3,14159**02423789 | **3,1415926**266579 |
| $n$ = 11 | **3,141**7360992607 | **3,14159**41599212 | **3,1415926**683944 |
| $n$ = 20 | **3,1415**657346587 | **3,1415925**761871 | **3,141592653**2636 |
| $n$ = 21 | **3,141**6160719183 | **3,1415927**142891 | **3,141592653**8114 |
| $n$ = 40 | **3,1415**890289487 | **3,14159265**11543 | **3,14159265358**7 |
| $n$ = 41 | **3,14159**60255683 | **3,14159265**57431 | **3,1415926535**923 |
| $n$ = 70 | **3,14159**19552651 | **3,141592653**4413 | **3,1415926535897** |
| $n$ = 71 | **3,14159**33232242 | **3,141592653**7284 | **3,141592653589**8 |

L'approximation du rapport de la circonférence au diamètre attribuée à Mādhava par ses disciples a certainement été obtenue en utilisant l'une de ces suites, mais il est aussi possible que les suites formées par les moyennes de deux termes consécutifs de chacune de ces suites aient été considérées afin de déterminer cette approximation plus rapidement encore. On constate que l'accélération de convergence est considérable, surtout eu égard à l'extrême lenteur de la convergence tant de la série initiale que du développement en fraction continue généralisée du reste. En ce qui concerne les suites $(b_n)$ et $(c_n)$, on peut même montrer que cette accélération est légèrement meilleure encore que celle induite par l'application réitérée de l'algorithme du « delta-2 d'Aitken ».

La mise en œuvre de la méthode exposée dans cet article serait-elle couronnée d'un succès comparable dans l'étude d'autres séries que celle de Mādhava-Leibniz ?

## Bibliographie